\documentclass{amsart}

\usepackage{amsthm}
\usepackage{amsmath}
\usepackage{amssymb}
\usepackage{amsfonts}
\usepackage{amscd}
\usepackage{eufrak}
\usepackage{textcomp}

\newtheorem{thm}{Theorem}
\newtheorem{cor}{Corollary}

\newtheorem*{rem}{Remark}
\newtheorem*{rems}{Remarks}
\newtheorem*{defi}{Definition}
\newtheorem{lem}{Lemma}

\begin{document}

	\title{On commuting exponentials in low dimensions}
	\author{GERALD BOURGEOIS}
	\address{Departement de Mathematiques, Faculte de Luminy, Route L.Lachamp, 13009 Marseille, France.}
	\email{bourgeoi@lumimath.univ-mrs.fr}
	\begin{abstract}
		$f, g \in \mathcal{L}(E)$ where $E$ is a $\Bbbk$ vector space of dimension $d$. We introduce the relation (*): 
		$\mbox{exp}(t.f+g)=\mbox{exp}(t.f) \circ \mbox{exp}(g)$
		 for any $t \in \Bbbk$ or $\mathbb N$; we study the connections between the relations ($^*$) and $f \circ g = g \circ f$ for $d=2$ or $3$. Let $d=2$: if $\Bbbk = \mathbb R$ and if (*) is verified for $t \in \mathbb N$ then $f \circ g = g \circ f$; we obtain all the couples $f$, $g$ verifying (*) on $\mathbb C$ and such as $f \circ g \neq g \circ f$.
		 Our main result is: if $d=3$, $\Bbbk = \mathbb C$ and 
		 $\mbox{exp}(t.f+g)=\mbox{exp}(t.f) \circ \mbox{exp}(g)=\mbox{exp}(g) \circ \mbox{exp}(t.f)$
		 for $t \in \mathbb N$, then $f$ and $g$ are simultaneously trigonalizable.
	\end{abstract}
	\maketitle

	\section{Introduction}
		Let $E$ be a vector space of dimension $d$ on $\Bbbk$ ( $\mathbb R$ or $\mathbb C$) and let $f,g \in \mathcal{L}(E)$.
		It is well known that $\{\forall t \in \Bbbk, \mbox{ exp}(t.(f+g))=\mbox{ exp}(t.f)\circ \mbox{ exp}(t.g)\}  \Rightarrow \{f \circ g = g \circ f\}$.
		However $A = 60 i \pi
		\left[
		\begin{tabular}{cc}
			$1$ & $0$ \\
			$0$ & $-1$
		\end{tabular}
		\right]$ 
		and 
		$B=  \pi 
		\left[
		\begin{tabular}{cc}
			$-150i$ & $-91$\\
			$391$ & $150i$
		\end{tabular}
		\right]$
	 	are not simultaneously trigonalizable and verify the following property: 
	 	$\forall n \in \mathbb N \mbox{ exp}(n.(A+B))=\mbox{exp}(n.A)\mbox{ exp}(n.B)$;
 		there is no contradiction with the preceding result because $\mathbb N$ has no accumulation point in $\Bbbk$.
 		
		In this paper we wish to show that the situation is quite different if one considers the relation(*): 
		$\mbox{ exp}(t.f+g)=\mbox{ exp}(t.f)  \circ \mbox{exp}(g) \mbox{ for } t \in \Bbbk \mbox{ or } \mathbb N$.

		\begin{rem} $A$ and $B$ verify (*) for $t=1,2,3,4,5$ but not for $t=6$.\end{rem}

		Several papers have dealt with $f,g$ such as
		\begin{enumerate}
			\item $\mbox{exp}(f+g)= \mbox{exp}(f) \circ \mbox{ exp}(g)$, with $f \circ g \neq g \circ f$
		\end{enumerate}
	
		Most remarkable are those of Morinaga and Nono: in \cite{2} the case where $d = 2$ is completely solved; in \cite{3} the case where $\Bbbk = \mathbb C$, $d=3$ is solved if 
		$\mbox{exp}(f) \circ \mbox{ exp}(g) = \mbox{exp}(g) \circ \mbox{ exp}(f)$ or if $f$ and $g$
		are simultaneously trigonalizable; these 2 types cover all the solutions for $d=2$ (see proof of theorem 2). Apparently M. \& N. have not tried to solve completely the problem for $d = 3$. We ignore if this question has been solved nowadays.
		
		If $d>3$ we do not know answers to the problem raised by (1), as well as to the following ones: what about $f,g$ if (2) or (3):
		\begin{enumerate}
		\value{enumi}1
			\item $e^f = e^g$
			\item $e^f \circ e^g = e^g \circ e^f$
		\end{enumerate}
		One unblocks the situation by using this notion:
		
		\begin{defi}
			$\mathcal{F} \subset \mathbb C$ is said to be 2i$\pi$\textit{-congruence free} (or incongruent mod 2i$\pi$) \linebreak
		\textit{iff} \mbox{ }$\forall u,v\in \mathcal{F}, u-v\notin 2i \pi \mathbb Z^*$.
		\end{defi}
		
		\begin{rem}
			the spectrum of f is 2i$\pi$-congruence free $\Leftrightarrow$ $f$ is a polynomial in exp(f).
		\end{rem}
		
		Here are the main results obtained with the above assumption:
		
		In \cite[p.161, lemma 8]{3}, M. \& N. show that if the spectrum of $f$ is 2i\textit{$\pi$-congruence free} and if $\mbox{exp}(f) = \mbox{ exp}(g)$, then $f \circ g = g \circ f$. Hille \cite{1} extends this result to the case of the bounded linear operators on a Hilbert space.
		
		In \cite[p.161, lemma 7]{3}, M. \& N. show that if the spectra of $f$ and $g$ are 2i$\pi$\textit{-congruence free} then $e^f \circ e^g = e^g \circ e^f \Rightarrow f \circ g = g \circ f$; Wermuth \cite{7} who did not know \cite{3} gives an alternative proof; Wermuth \cite{8} extended this result to the situation where $f,g$ are bounded linear operators on a Banach space. Schmoeger \cite{5} simplified this last proof.
		
		The equation $e^f \circ e^g = e^{f+g}$ is more difficult to handle except if $f$ and $g$ are bounded self-adjoint operators on a Hilbert space. It is an important test of commutativity, nevertheless, as shown by this result:
		
		if $\mathcal{A}$ is a $C^{\ast}$-algebra, then: $\mathcal{A}$ is commutative $\Leftrightarrow \forall f, g$ positive $\in \mathcal{A}$, $e^f \circ e^g = e^{f+g}$ in $\tilde{\mathcal{A}}$ (\textit{cf.} \cite{9}).
		
		If the spectrum of $f$ is 2i\textit{$\pi$-congruence free}, Schmoeger \cite{6} showed, within the framework of the bounded operators, that if $e^f \circ e^g = e^g \circ e^f$, then $f \circ g - g \circ f$ is a sum of nilpotents; if moreover $e^f \circ e^g = e^{f+g}$ then $f \circ g = g \circ f$ (Paliogiannis \cite{4} gives an alternative proof).
		
		~
		
		Thus there are 2 ways of broaching the relations (1),(2),(3):
		
		\begin{enumerate}
		\renewcommand{\labelenumi}{\alph{enumi})}
			\item One uses the 2i\textit{$\pi$-congruence free} hypothesis: the advantage of this hypothesis is that it provides results which are valid even for bounded linear operators. It has a serious drawback, however: this way of calculus rejects \textit{a priori} a large part of the $f,g$ verifying the relation under study. Currently we observe quite a consensus in favour of this point of view.
			\item On the contrary if we don't use this hypothesis then we can deal with all the $f,g$ pairs. But now one is restricted to $d \leq 3$. This was M. \& N.'s position although they knew the power of the hypothesis over the spectra.
		\end{enumerate}		
			
		In this paper we adopt the second position. We draw on \cite{2} and \cite{3}, in order to study the relation($^*$). 
		We use the particular case $d = 2$ to prove theorem 3 ( where $d = 3$), which constitutes our main result.
		Let us notice that this last theorem would have a trivial conclusion in the conditions of a).

	\section{Dimension 2}
	
	\begin{rems}
	
	$ $
	
		\begin{enumerate}
		\renewcommand{\labelenumi}{\alph{enumi})}
			\item $id$ refers to identity on $E$. If $u \in \mathcal{L}(E)$, $tr(u)$ stands for the trace of u.
			\item We will use repeatedly this piece of calculus:
		\end{enumerate}
	\end{rems}
	
	\begin{lem}
		 Let $P(T)\in \mathbb Z[T]$ the polynomial $P(T) = \alpha^2 T^2 + \beta T + \gamma$ where $\alpha \in \mathbb N^*$.
		$(t_n)_n$ is a stricly increasing series of integers.
		
		If \mbox{ }$\forall n \in \mathbb N$ $P(t_n)$ is a square then $P$ is the square of a one degree polynomial \textit{i.e.} $\beta^2 = 4\alpha^2\gamma$.
	\end{lem}
	
	\begin{proof}
		$\forall n \in \mathbb N^*$ $\exists u_n \in \mathbb N^*$ such as 
		$P(t_n)={u_n}^2=\alpha^2 {t_n}^2 (1+ \frac \beta {\alpha^2 t_n} + \frac {\gamma}{\alpha^2{t_n}^2})$;
		\newline
		when $t_n \rightarrow \infty$ $u_n = \alpha t_n + \frac {\beta} {2\alpha}+O(\frac{1}{t_n})$, which implies that:
		\newline
		$\frac{\beta}{2\alpha} \in \mathbb Z$ and $u_n = \alpha t_n + \frac{\beta}{2\alpha}$ for a quite large $n$. It follows $\frac{\beta^2}{4\alpha^2}=\gamma$
	\end{proof}
	
	\begin{thm}
		\label{thm3}
			If $\Bbbk = \mathbb R$ and $d = 2$, let a strictly increasing series of integers $(t_n)_n$ such as $t_0=1$; then 
			$\{\forall n \in \mathbb N, \mbox{exp}(t_n f+g) = \mbox{exp}(t_nf) \circ \mbox{exp}(g)\} \Leftrightarrow \{f \circ g = g \circ f\}$.
	\end{thm}
	
	\begin{proof} \textit{N.B.} :
	~
		\begin{enumerate}
		\renewcommand{\labelenumi}{\alph{enumi})}
			\item only ($\Rightarrow$) has to be shown.
			\item  if (*) holds for $f,g$ then $\forall \sigma, \tau \in \mathbb C$ (*) holds for ($f-\sigma.id, g -\tau.id$).
		\end{enumerate}
		From now on we suppose that if ($f,g$) verify (*) then $tr(f) = tr(g) = 0$.
		\newline
		Here we reason \textit{ad absurdum}. We first recall, in our terms, one of M. \& N's results (\cite[p. 357]{2}):
		
		$\{e^{f+g}=e^f \circ e^g, f \circ g \neq g \circ f, tr(f)=tr(g)=0\}$
		\newline
		implies that there is a $\mathbb R^2$ basis in which $f$ ang $g$ have as representative matrices:
		$A = \pi
		\left[
		\begin{tabular}{cc}
			$0$ & $-\lambda$ \\
			$\lambda$ & $0$
		\end{tabular}
		\right]$ 
		where $\lambda \in \mathbb N^*$ and 
		$B = \pi
		\left[
		\begin{tabular}{cc}
			$a$ & $b$ \\
			$c$ & $-a$
		\end{tabular}
		\right]$
		 where $a,b,c \in \mathbb R$ such as:
		 
		\begin{enumerate}
		\value{enumi}3
			\item $spectrum(B) = \{i \pi \mu, -i \pi \mu\}$ and $spectrum(A+B)=\{i \pi \nu, -i \pi \nu\}$ where $\mu, \nu \in \mathbb N^*$.
			\item $\nu^2 \neq (\lambda \pm \mu)^2$
			\item $(b \neq -c$ or $a \neq 0)$ and $-a^2-bc=\mu^2$.
		\end{enumerate}
		
		\begin{rems}
		~
			\begin{enumerate}
			\renewcommand{\labelenumi}{\alph{enumi})}
				\item $e^A.e^B = e^B.e^A$
				\item The spectra of $t_n.A, B$ and $t_n.A+B$ are never 2i$\pi$-congruence free.
			\end{enumerate}
		\end{rems}
		Clearly $\nu^2=\frac{1}{\pi^2}$det$(A+B)=\lambda^2+\mu^2-\lambda(b-c)$.
		\newline
		The representative matrices of $t_nf$ and $f$ have the same form; this implies that
		$(t_n \lambda)^2+ \mu^2 -(t_n\lambda)(b-c)=\lambda^2 t_n^2 +t_n(\nu^2 - \lambda^2 - \mu^2)+ \mu^2$ is a square.
		
		From lemma 1 it follows that 
		$| \nu^2 - \lambda^2 - \mu^2|=2 \lambda \mu$ and $\nu^2=(\lambda \pm \mu)^2$ which is the contradictory of (5).	
	\end{proof}
	
	Let $\mathcal{U} = \{u \in \mathbb C^*|\mbox{ }e^u = 1+u\}$; this set contains an infinity of elements including $u \approx 2.0888 + 7.4615 i$.
	
	\begin{thm}
		\label{thm3}
			If $\Bbbk = \mathbb C$ and $d = 2$, then the following equivalence is true:
			\newline
			$\{ f \circ g \neq g \circ f \mbox{ and } \forall t \in \Bbbk$ exp$(t.f+g)=exp(t.f) \circ exp(g)\}$
			$\Leftrightarrow \{ \exists \sigma, \tau \in \mathbb C \mbox{ such as: } \tilde{f}=f-\sigma.id \mbox{ and } \tilde{g}=g-\tau.id \mbox{ verify }: \tilde{f}\neq 0, \tilde{f^2} = 0, \tilde{f} \circ \tilde{g}= 0 \mbox{ and } tr(\tilde{g}) \in \mathcal{U}\}$.
	\end{thm}
	
	\begin{proof} \textit{N.B.}: the content of $f,g$ is indifferent to the addition of an homothety.

	$(\Rightarrow)$: according to \cite[p.356]{2}, the study of ($f,g$) such as $f \circ g \neq g \circ f$ and exp($f+g$)= exp($f$)$\circ$ exp($g$) reduces to 4 cases:
	
	\textit{1$^\circ$ case}: there is a $\mathbb C^2$ basis in which the matrices have the same form as in the real case (with exp($f$)$\circ$ exp($g$) = exp($g$)$\circ$ exp($f$), \textit{cf} the proof of theorem 1), but here $a,b,c \in \mathbb C$; moreover one can have $\nu^2=(\lambda \pm \mu)^2$
		and then [$f,g$] is a linear combination of $f,g$ and there is a $\mathbb C^2$ basis in which $f$ and $g$ have the representative matrices:
		$A =
		\left[
		\begin{tabular}{cc}
			$i \pi \lambda$ & $0$ \\
			$0$ & $-i \pi \lambda$
		\end{tabular}
		\right]$ 
		and
		$B =
		\left[
		\begin{tabular}{cc}
			$i \pi \mu$ & $1$ \\
			$0$ & $-i \pi \mu$
		\end{tabular}
		\right]$
		where $\lambda, \mu \in \mathbb Z^*$ and $\lambda + \mu \neq 0$.
		
		The preceding proof shows that $\nu^2=(\lambda \pm \mu)^2$ is necessarily true; ($^*$) holds for the above couple ($A,B$) with $t \in \mathbb N$ \textit{but not} for $t \in \mathbb C$ such as $\lambda t \notin \mathbb Z$.
		
		~
		
		In the last 3 cases, $f$ and $g$ are simultaneously trigonalizable and $e^f \circ e^g \neq e^g \circ e^f$:
		
		~
		
	\textit{2$^\circ$ case}: there is a $E$ basis in which the representative matrices of $f$ and $g$ are 
		$A =
			\left[
			\begin{tabular}{cc}
				$0$ & $0$ \\
				$0$ & $u$
			\end{tabular}
			\right]$ 
		and 
		$B =
			\left[
			\begin{tabular}{cc}
				$v$ & $1$ \\
				$0$ & $0$
			\end{tabular}
			\right]$ 
		with $u \in \mathbb C^*$, $u \neq v$, and $(\{v \neq 0 \mbox{ and } \frac{e^u-1}{u} = \frac{e^v-1}{v} \neq 0\} \mbox{ or }\{v = 0 \mbox{ and } u \in \mathcal{U} \})$. 
		
		($^*$) holds for $t=2$ under the condition \{$v \neq 0$ and $(e^u-1)^2 v  = 0$\}$ \mbox{ or }\linebreak $\{$v = 0 \mbox{ and } u^2 = 0$\}, which is absurd.
		
	\textit{3$^\circ$ case}: there is a $E$ basis in which the representative matrices of $f$ and $g$ are
		$A =
			\left[
			\begin{tabular}{cc}
				$u$ & $0$ \\
				$0$ & $0$
			\end{tabular}
			\right]$ 
		and 
		$B =
			\left[
			\begin{tabular}{cc}
				$0$ & $1$ \\
				$0$ & $u$
			\end{tabular}
			\right]$ 
		with $u \in \mathcal{U}$.
		
	($^*$) holds for $t=2$ under the condition $e^u u = 0$, which is absurd.
	
	\textit{4$^\circ$ case}: there is a $E$ basis in which the representative matrices of $f$ and $g$ are 
		$A =
			\left[
			\begin{tabular}{cc}
				$0$ & $1$ \\
				$0$ & $0$
			\end{tabular}
			\right]$ 
		and 
		$B =
			\left[
			\begin{tabular}{cc}
				$u$ & $0$ \\
				$0$ & $0$
			\end{tabular}
			\right]$ 
		where $u \in \mathcal{U}$.
	
	$AB = 0$ and $e^{t.A} = id +t.A$ from \textit{spectrum}$(t.A + B)=\{u, 0\}$, we deduce $e^B = id + B$, $e^{t.A+B} = id + t.A + B$; this makes it possible to show that: $\forall t \in \mathbb C^*$ $e^{t.A}e^B = e^{t.A+B}(\neq e^Be^{t.A}$).
	
	This is the unique case which provides a solution; we verify that $A \neq 0$, $A^2 = 0$, $AB = 0$, $tr(B) = u$.
	
	($\Leftarrow$): Conversely if one starts from these 4 relations, one is driven by a first basis change to: $A =
			\left[
			\begin{tabular}{cc}
				$0$ & $1$ \\
				$0$ & $0$
			\end{tabular}
			\right]$ 
		and 
		$B =
			\left[
			\begin{tabular}{cc}
				$u$ & $\theta$ \\
				$0$ & $0$
			\end{tabular}
			\right]$ , 
	then by a second basis change to the 2 type-matrices of the 4$^\circ$ case; then ($^*$) holds for any t.
	\end{proof}		
	\begin{rem}
		the spectra of t.f, g and t.f+g are always $2i\pi$-congruence free.
	\end{rem}
	
	$ $
	
	\begin{cor}
		Consider the same framework $\Bbbk = \mathbb C$ and $d = 2$. Now suppose a strictly increasing series of integers $(t_n)_n$ such as $t_0 = 1$ and $t_1 = 2$; it follows that, if $ \forall n \in \mathbb N$, $exp(t_nf + g) = exp(t_nf) \circ exp(g)$, then $f$ and $g$ are simultaneously trigonalizable.
	\end{cor}
	
	\begin{proof}
		If $f \circ g = g \circ f$ then $f$ and $g$ are simultaneously trigonalizable.

		Assume $f \circ g \neq g \circ f$; then we follow the proof of theorem 2: only the 1$^\circ$ and 4$^\circ$ cases provide solutions; in each one of these cases $f$ and $g$ are simultaneously trigonalizable.
	\end{proof}

	\section{Dimension 3}

		\begin{thm}
			\label{thm3}
			Here $\Bbbk = \mathbb{C}$ and $d=3$; suppose a strictly increasing series of integers $(t_n)_n$ such that $t_0=1$ and $t_1=2$; if \mbox{ }$\forall n \in \mathbb{N}$, $\exp(t_n.f+g)=\exp(t_n.f) \circ \exp(g) = \exp(g) \circ \exp(t_n.f)$, then $f$ and $g$ are simultaneously trigonalizable.
		\end{thm}
		
		$ $
		
		\begin{rems}
			~
			\begin{enumerate}
			\renewcommand{\labelenumi}{\alph{enumi})}
				\item Schmoeger's already mentioned result \cite{6} implies that, if $f$, $g$ are defined as in theorem \ref{thm3}, then they commute or else their spectra are not $2i\pi$-congruence free.
				\item In the proof we suppose $t_n=n$ for the sake of simplicity.
			\end{enumerate}
		\end{rems}
		
		\begin{proof}
	If $u \in \mathcal{L}(E)$, we note \#($u$) the number of distinct eigenvalues of $u$. What follows makes it possible to determine all the ($f$, $g$) verifying the required conditions.
	If \#$(f)=3$, $\exists n_0 \in N$ such as $n \geq n_0 \Rightarrow \#(n.f) = \#(n.f + g) = 3$; even if it means to replace $f$ by $n_0.f$, we can then suppose 
	
	$(P): \{ \mbox{or } \#(f) \in \{1,2\}$ or $\forall n \in \mathbb N^*$ $\#(n.f) = \#(n.f + g) = 3\}$.
	
	In \cite[p. 164-177]{3}, the solutions of $e^f \circ e^g = e^g \circ e^f = e^{f+g}$ fall into 9 types; 2 types do not respect $(P)$; 3 other types provide simultaneously trigonalizable couples.
	
	There remain 4 types to be examined:
	
	1$^\circ$ \textit{case} (\textit{cf.} \cite[p. 175, theorem 7, case I]{3}): there is $n_0 \in \mathbb N^*$ and a $\mathbb C^3$ basis in which $n_0.f$ and $g$ are represented by 
	$A = \left[
		\begin{tabular}{cc}
			$X$ & $0$ \\
			$0$ & $\lambda$
		\end{tabular}
		\right]$ 
	and $B = \left[
		\begin{tabular}{cc}
			$Y$ & $0$ \\
			$0$ & $\mu$
		\end{tabular}
		\right]$ 
	where $X,Y \in \mathfrak{M}_2(C)$.
	
	$X,Y$ verify the assumptions of corollary 1 and thus $A,B$ are simultaneously trigonalizable.
	
	~
	
	In the last 3 cases exp($f$) = exp($g$) = exp($f+g$) = $id$. We overlook, of course, the ($f$, $g$), belonging to the 1$^\circ$ case.
	
	~
	
	2$^\circ$ \textit{case} (\textit{cf.} \cite[p. 173-175, case III$_4$]{3} $\forall n \in \mathbb N^*$ $\#(f) = \#(f+g) = 3$.
	
	There is a $\mathbb C^3$ basis in which $\frac{1}{2i\pi}f$ and $\frac{1}{2i\pi}g$ are represented by 
	
	$A = s^{-1}.diag(l_1,l_2,0).s$, $ B = diag(m_1,m_2,m_3)$, such as 
	
	$ A + B = t^{-1}.diag(n_1,n_2,0).t$ where
	
	$l_1,l_2 \in \mathbb Z^*, l_1\neq l_2, m_1,m_2,m_3 \in \mathbb Z, m_1 \neq m_2, n_1,n_2 \in \mathbb Z^*, n_1 \neq n_2 $.
	
	If $A = [a_{ij}]$, \cite[p. 174]{3} indicates: $ \exists \mbox{ }\rho,\sigma \in \mathbb C $ such as:

	\begin{enumerate}
		\value{enumi}6
			\item $a_{33} = \rho (m_1 - m_2)$,
			\item $a_{12} = \rho (m_1^2 - m_2^2) + \sigma ( m_1 - m_2 )$,
			\item $a_{23} = \rho (m_2^2 - m_3^2) + \sigma ( m_2 - m_3 ) - \frac{1}{m_1 - m_2}[(m_2 - m_3)l_1l_2 + m_3(m_3 - n_1)(m_3 - n_2)]$,
			\item $a_{31} = \rho (m_3^2 - m_1^2) + \sigma ( m_3 - m_1 ) + \frac{1}{m_1 - m_2}[(m_1 - m_3)l_1l_2 + m_3(m_3 - n_1)(m_3 - n_2)]$,
			\item $a_{11} = \rho (m_2 - m_3) + \frac{1}{m_1 - m_2}[(l_1 + l_2)(m_1+m_3)+l_1l_2+(m_1m_2+m_2m_3+m_3m_1)-n_1n_2]$,
			\item $a_{22} = \rho (m_3 - m_1) - \frac{1}{m_1 - m_2}[(l_1 + l_2)(m_2+m_3)+l_1l_2+(m_1m_2+m_2m_3+m_3m_1)-n_1n_2]$,
		\end{enumerate}

	If one changes $A$ into $nA$, the equalities are preserved by changing $a_{ij}$ into $n.a_{ij}$, $l_i$ into $n.l_i$, $m_i$ into $m_i+\lambda$, $n_i$ into $\tilde{n_i}$ ( we choose $\lambda \in \mathbb Z$ so that $\tilde{n_3}=0$), $\rho$ into $\tilde{\rho}$, $\sigma$ into $\tilde{\sigma}$.
		
	By (7)\mbox{ } $\tilde{\rho}=n\rho$; by (8)\mbox{ } $\tilde{\sigma}=n(\sigma - 2 \lambda \rho)$.
	
	The other 4 equations imply $n=0$ or $n=1$ and thus ($f$, $g$) is inappropriate; this point follows from a calculus using the Maple "Grobner" package.
	
	~
		
	3$^\circ$ \textit{case} (\textit{cf.} \cite[p. 171, case III$_2$(i)]{3}): \#($f$) = 2, \#($g$) = \#($f+g$) = 3.
	
	There is a basis in which $\frac{1}{2i\pi}f$ and $\frac{1}{2i\pi}g$ are represented by 
	
	$A=[a_{ij}]= s^{-1}.diag(l_1,0,0).s$, $B = diag(m_1,m_2,m_3)$ such as: 
	
	$A+B = t^{-1}.diag(n_1,n_2,0).t$ \mbox{ } where:
	
	$l_1 \in \mathbb Z^*$; $m_1, m _2,m_ 3 \in \mathbb Z$ and are distincts 2 by 2; 
	
	$n_1,n_2 \in \mathbb Z^*$ and $n_1 \neq n_2; m_1+m_2+m_3 \neq n_1+n_2$;
	
	$a_{11}= \frac{m_1(m_1-n_1)(m_1-n_2)}{(m_1-m_2)(m_3-m_1)}$, $a_{22} = \frac{m_2(m_2-n_1)(m_2-n_2)}{(m_2-m_3)(m_1-m_2)}$, $a_{33} = \frac{m_3(m_3-n_1)(m_3-n_2)}{(m_3-m_1)(m_2-m_3)}$.
	
	\begin{rem}
		if \mbox{ }$\exists n \in \mathbb N^*$ such that $\#(n.f+g) \neq 3$, then ($f$, $g$) does not belong to any of the 4 cases and is inappropriate; one can then suppose that $\forall n \in \mathbb N^* (n.f,n.f+g)$ falls into the 3$^\circ$ case.
	\end{rem}
	
	Now we inspect 2 subcases:
	
	i) $a_{11}a_{22}a_{33} \neq 0$; then one returns to $A = \left[
		\begin{tabular}{ccc}
			$a_{11}$ & $\sqrt{a_{11}a_{22}}$ & $\sqrt{a_{11}a_{33}}$ \\
			$\sqrt{a_{11}a_{22}}$ & $a_{22}$ & $\sqrt{a_{22}a_{33}}$ \\
			$\sqrt{a_{11}a_{33}}$ & $\sqrt{a_{22}a_{33}}$ & $a_{33}$			
		\end{tabular}
		\right]$ 
		where the square roots are selected such as rank $(A)=1$; $\forall n \in \mathbb N^*$:
		
		$det(nA + B - x.id)=x^3 + ((n-1)(m_1+m_2+m_3)-n(n_1+n_2))x^2+\linebreak(-(n-1)(m_1m_2+m_2m_3+m_3m_1)+nn_1n_2)x+(n-1)m_1m_2m_3$
		
		has its roots $u_n,v_n,w_n$ in $\mathbb Z^*$ and such as: $|u_n| \geq |v_n| \geq |w_n| \geq 1$ ( because $ m_1m_2m_3 \neq 0 $).
		
		If $\sigma_1,\sigma_2,\sigma_3$ are their symmetrical functions:
		
		$\sigma_1 \sim n(n_1+n_2-m_1-m_2-m_3)$; $\sigma_2 = O(n)$; $\sigma_3 \sim -nm_1m_2m_3$; ${u_n}^2+{v_n}^2+{w_n}^2 \sim {\sigma_1}^2 = \Theta(n^2)$ thus $u_n = \Theta(n)$ and $v_n,w_n$ are bounded;

		$u_n \sim \sigma_1$, $v_nw_n \rightarrow \lambda = \frac{m_1m_2m_3} {m_1+m_2+m_3-n_1-n_2}$, $v_n+w_n \rightarrow \mu = \frac{m_1m_2+m_2m_3+m_3m_1-n_1n_2} {m_1+m_2+m_3-n_1-n_2}$.
		
		As we work in $\mathbb Z$, for a quite large $n: v_nw_n = \lambda$, $v_n + w_n = \mu$.
		\begin{itemize}
			\item If $\mu \neq 0$
			\newline
			$u_n = \sigma_1 - \mu = \frac{\sigma_3}{\lambda} = \frac{\sigma_2-\lambda}{\mu}$ hence $\mu = n_1 + n_2, \lambda = n_1n_2$; then for example $\lambda = m_2m_3$, $\mu = m_2 + m_3$; it results from it that 
			$(m_2 - n_1)(m_2 - n_2) = 0$ and $a_{22} = 0$ which is a contradiction.
			
			\item If $\mu = 0$
			\newline
			$v_n+w_n = 0$, $u_n = \sigma_1 = \frac{\sigma_3} {\lambda}$ from where $n_1+n_2 = 0$; $\sigma_2 = v_nw_n = \lambda$ is constant thus its value is $m_1m_2 + m_2m_3 + m_3m_1$.		
			\newline
			Then $\exists k \in \mathbb Z^*$ such as $m_1m_2 +m_2m_3 +m_3m_1 = \frac{m_1m_2m_3} {m_1+m_2+m_3} = -k^2$, from where for example $k = m_1 = -m_2$.
			\newline
			Then $n_1n_2 = m_1m_2 + m_2m_3 + m_3m_1 = m_1m_2$; thus $(m_1 - n_1)(m_1 - n_2) = 0$ which is contradictory.
			\end{itemize}

		ii) $a_{11}a_{22}a_{33} = 0$; for example $a_{33} = 0$; there are 8 possible forms for $A$, 2 of which are triangular and 2 return to the 1$^\circ$ case; 4 forms $(A_i)_{1 \leq i \leq 4}$ remain to be examined:

		$A_1 = \left[
		\begin{tabular}{ccc}
			$a_{11}$ & $\sqrt{a_{11}a_{22}}$ & $\sqrt{a_{11}a_{22}}$\\
			$\sqrt{a_{11}a_{22}}$ & $a_{22}$ & $a_{22}$\\
			$0$ & $0$ & $0$
		\end{tabular}
		\right]$,
		where $a_{11} = \frac{(m_1-n_1)(m_1-n_2)} {m_2-m_1} \neq 0$, $a_{22} = \frac{(m_2-n_1)(m_2-n_2)}{m_1-m_2} \neq 0$; $m_3$ is an eigenvalue of $n.A + B$; one can thus suppose that $m_3 = 0$ and $m_1 +m_2 \neq n_1 +n_2$; the other 2 eigenvalue are the roots (in $\mathbb Z$) of:
		
		$x^2 + ((n - 1)(m_1 + m_2) - n(n_1 + n_2))x + (1 - n)m_1m_2 + nn_1n_2$;
		
		$\forall n \in \mathbb N^*$ the discriminant:		
		
		${(m_1+m_2-n_1-n_2)}^2n^2+2((m_1+m_2)(n_1+n_2)-{m_1}^2-{m_2}^2-2n_1n_2)n+{(m_1-m_2)}^2$ must be a square.
		
		By lemma 1 (as $m_1 + m_2 \neq n_1 + n_2$) this polynomial in $n$ is the square of one degree polynomial and its discriminant is null: 
		
		$(m_1 -n_1)(m_1 -n_2)(m_2 -n_1)(m_2 - n_2) = 0$, which is contradictory.
		
		$A_2 = {A_1}^*$ does not hold by the same calculus.
		
		$A_3 = \left[
		\begin{tabular}{ccc}
			$l_1$ & $1$ & $0$\\
			$0$ & $0$ & $0$\\
			$l_1$ & $1$ & $0$
		\end{tabular}
		\right]$
		and $B$ are simultaneously trigonalizable;
		
		the same result holds for 
		$A_4 = \left[
		\begin{tabular}{ccc}
			$l_1$ & $0$ & $1$\\
			$l_1$ & $0$ & $1$\\
			$0$ & $0$ & $0$
		\end{tabular}
		\right]$
		and $B$.
		
		\newpage
		
		4$^\circ$ \textit{case} (\textit{cf.} \cite[p.172, case III$_2$(ii)]{3}): $\#(f) = 2,\#(g) = 2,\#(f + g) = 3$.

		There is a basis in which $\frac{1}{ 2i\pi }f$ and $\frac{1}{ 2i\pi }g$ are represented by: 
		
		$A = [a_ib_j ] = s^{-1}.diag(l_1, 0, 0).s$, $B = diag(m, 0, 0)$ such as 
		
		$A + B = t^{-1}.diag(n_1, n_2, 0).t$ where 
		
		$l_1,m, n_1, n_2 \in \mathbb Z^*$, $n_1 \neq n_2$, $m \neq n_1 + n_2$ and finally: $\exists \alpha \in \mathbb C$ such as 
		
		$a_1b_1 = \frac{-(m-n_1)(m-n_2)}{m}$ , $a_2b_2 = -\alpha + \frac {n_1n_2}{m}$ , $a_3b_3 = \alpha$.
		
		Here also we can suppose that $\forall n \in  \mathbb N^*$ $(n.f, n.f + g)$ meets the conditions of the 4$^\circ$ case.
		
		$\forall n \in  \mathbb N^*$ $nA + B$ has a null eigenvalue, the other two being roots (in $\mathbb Z$) of:
	
		$x^2 - (n(a_1b_1 + a_2b_2 + a_3b_3) + m)x + mn(a_2b_2 + a_3b_3)$;
		
		$\forall n \in  \mathbb N^*$ the discriminant $(a_1b_1+a_2b_2+a_3b_3)^2n^2+2m(a_1b_1-a_2b_2-a_3b_3)n+m^2$ must be a square;
		
		by lemma 1 (as $a_1b_1+a_2b_2+a_3b_3 = n_1+n_2-m \neq 0$), this polynomial in $n$ is the square of one degree polynomial and its discriminant is null:

		$(a_2b_2 + a_3b_3)m^2a_1b_1 = n_1n_2ma_1b_1 = 0$.
		
		For example $a_1 = 0$ from where rank $(AB - BA) = 1$ and $A$ and $B$ are simultaneously trigonalizable.
\end{proof}
		
	\section{Conclusion}

	When one forces the spectrum of $f$ to be $2i\pi$-\textit{congruence free}, then $z \rightarrow e^z$ is
one to one on a neighborhood $D$ of the spectrum of $f$ and one can (as in \cite{4}) make
a pure logical reasoning using of the holomorphic functions on $D$. In theorem 3 we
do not make this assumption on $f$ and we reason by examinations of cases. This
method cannot work in dimension $\geq 4$; can one still make a pure logical reasoning
if the exponential is not invertible?

We may only hope that, in the future, it will be possible to demonstrate our
\textit{conjecture}:

Theorem 3 is valid in \textit{any} dimension.

\end{document}